\documentclass[a4paper]{amsart}

\author{Samuel Boissi{\`e}re}

\title{Chern classes of the tangent bundle on the Hilbert scheme of
points on the affine plane}

\address{Samuel Boissi\`{e}re, Fachbereich f\"{u}r Mathematik, Staudinger Weg $9$,
Johannes Gutenberg-Universit\"{a}t Mainz, $55099$ Mainz, Deutschland}

\email{boissiere@mathematik.uni-mainz.de}

\urladdr{http://sokrates.mathematik.uni.mainz.de/$\sim$samuel}

\usepackage{amsmath,amssymb}
\usepackage[all]{xy}
\usepackage[english]{babel}
\usepackage{eufrak,mathrsfs}
\usepackage{youngtab}

\DeclareMathOperator{\Hilb}{Hilb} \DeclareMathOperator{\Coeff}{Coeff}
\DeclareMathOperator{\vac}{|0\rangle} \DeclareMathOperator{\End}{End}
\DeclareMathOperator{\Supp}{Supp} \DeclareMathOperator{\id}{id}
\DeclareMathOperator{\ad}{ad} \DeclareMathOperator{\rk}{rk}
\DeclareMathOperator{\Bl}{Bl}

\newcommand{\IC}{\mathbb{C}}
\newcommand{\IH}{\mathbb{H}}
\newcommand{\IN}{\mathbb{N}}
\newcommand{\IQ}{\mathbb{Q}}
\newcommand{\IZ}{\mathbb{Z}}

\newcommand{\hilb}{\Hilb^n(\IC^2)}

\newcommand{\cC}{\mathcal{C}}

\newcommand{\cI}{\mathcal{I}}
\newcommand{\cL}{\mathcal{L}}
\newcommand{\cO}{\mathcal{O}}
\newcommand{\cU}{\mathcal{U}}

\newcommand{\cD}{\mathscr{D}}

\newcommand{\kc}{\mathfrak{c}}
\newcommand{\kch}{\mathfrak{ch}}
\newcommand{\kd}{\mathfrak{d}}
\newcommand{\kf}{\mathfrak{f}}
\newcommand{\kg}{\mathfrak{g}}
\newcommand{\kq}{\mathfrak{q}}

\newcommand{\kF}{\mathfrak{F}}
\newcommand{\kG}{\mathfrak{G}}
\newcommand{\kL}{\mathfrak{L}}
\newcommand{\kU}{\mathfrak{U}}

\newcommand{\ie}{\textit{i.e. }}

\theoremstyle{plain}
\newtheorem{theorem}{Theorem}[section]
\newtheorem{definition}[theorem]{Definition}
\newtheorem{proposition}[theorem]{Proposition}
\newtheorem{lemma}[theorem]{Lemma}
\newtheorem{remark}[theorem]{Remark}
\newtheorem{corollary}[theorem]{Corollary}

\begin{document}

\begin{abstract}
The cohomology of the Hilbert schemes of points on smooth projective surfaces can be approached both with vertex algebra tools and equivariant tools. Using the first tool, we study the existence and the structure of universal formulas for the Chern classes of the tangent bundle over the Hilbert scheme of points on a projective surface. The second tool leads then to nice generating formulas in the particular case of the Hilbert scheme of points on the affine plane.
\end{abstract}



\maketitle

\pagestyle{myheadings}

\markboth{SAMUEL BOISSI\`{E}RE}{Chern classes of the tangent bundle on the Hilbert scheme of points
on the affine plane}

\section{Introduction}

Let $S$ be a smooth projective complex surface and $n$ a non-negative integer. The Hilbert scheme of $n$ points on $S$, denoted by $\Hilb^n(S)$ or $S^{[n]}$, parameterizing generalized $n$-tuples on $S$, \ie zero-dimensional subschemes of length $n$ of $S$ is, by a result of Fogarty \cite{Fo}, a smooth projective surface of complex dimension $2n$. In the study of his rational cohomology, initiated by G{\"ottsche} \cite{Goe}, Grojnowski \cite{Gro} and Nakajima \cite{N3}, several vertex algebra tools have been developed by Lehn \cite{L} and after by Li-Qin-Wang \cite{LQW3,LQW2,LQW4}, providing a better understanding of the universality of the ring structure of the cohomology as such as universal formulas for the Chern classes of tautological bundles on the Hilbert scheme.

In a different flavor, the equivariant structure of the Hilbert scheme of points on the projective plane for a natural action of the torus has been used by Ellingsrud-Str{\o}mme \cite{ES}, Nakajima \cite{N2} and Vasserot \cite{V} to get a better understanding of the tangent bundle and the cohomology ring  of $\hilb$. In this context, the  equivariant cohomology of the Hilbert scheme of points on the affine plane was the main tool in our preceding work \cite{SB1} to get general combinatorial formulas for the Chern classes of linearized vector bundles over $\hilb$. 

In view of the existing generating formulas for the Chern classes of tautological bundles discovered by Lehn \cite{L}, it is natural (as we learned with Lehn \cite{L2}) to ask for a similar result in the case of the tangent bundle over the Hilbert scheme of points on the affine plane. The total cohomology space $\bigoplus\limits_n H^*(\hilb)$ is naturally isomorphic, as a vertex algebra, to the space of polynomials $\Lambda:=\IQ[p_1,p_2,p_3,\ldots]$. In this paper, we prove the following generating formulas for the Chern classes and the Chern characters of the tangent bundles $T\hilb$, expressed in $\Lambda$:

\begin{theorem}
The Chern classes of the tangent bundle over the Hilbert scheme $\hilb$ are given
by the following generating series:
$$
\sum_{n\geq 0} c\left(T\hilb\right)=\exp\left( \sum_{k\geq 0} (-1)^k C_k
\frac{p_{2k+1}}{2k+1} \right),
$$
where $C_k:=\displaystyle \frac{1}{k+1}\binom{2k}{k}$ is the $k$-th Catalan
number.

The Chern characters are given by the following generating series:
$$
\sum_{n\geq 0} ch\left(T\hilb\right)= 2e^{\displaystyle p_1}\sum_{k\geq
0}\displaystyle\frac{p_{2k+1}}{(2k+1)!}.
$$
\end{theorem}

The strategy of the proof is as follows. We first
have to understand the common features of such generating series of Chern classes for the tangent bundle over $S^{[n]}$. To do so, we develop the notion of ``universal formula'' in the same spirit as in \cite{LQW2} and we prove a raffined version of the theorem of Ellingsrud-G{\"o}ttsche-Lehn \cite{EGL} on the Chern numbers of the tangent bundles $TS^{[n]}$ by establishing the universality of the Chern classes (\S \ref{ss:univTang}). Applied properly to the case of the affine plane, we see that the formulas are then much easier and we finish the computation by means of our general formulas \cite{SB1} and some combinatorial tricks.

This research is part of my thesis \cite{SB2} written at the University of Nantes. I owe special thanks to my supervisor C. Sorger for his continuous encouragement and the frequent discussions we had during these three years. I am very grateful to M. Lehn for all the things I learned from conversations with him, for his interest and support. I thank M. Nieper-Wi\ss kirchen for a very instructive correspondence on hypergeometric functions, and the anonymous referee for his carefully reading and his suggestions.

\section{Cohomology of Hilbert schemes of points}

In this section, we recall some classical results related to the vertex
algebra structure of the total cohomology space of Hilbert schemes of points on
surfaces.

{\it Conventions.} For a smooth quasi-projective variety $X$, we denote by $H^*(X)$ the singular cohomology with \emph{rational} coefficients of the underlying real manifold, by $K(X)$ the \emph{rational} Grothendieck group generated by locally free sheaves, or equivalently by arbitrary coherent sheaves and by $ch:K(X)\xrightarrow{\sim} H^{2*}(X)$ the Chern character. Our conventions and notations for the various operations in both theories follow \cite{BS}. For any continuous map $f:X\rightarrow Y$ between two smooth oriented compact manifolds, we define a cohomological push-forward by $f_!:=D_Y^{-1}\circ f_*\circ D_X$ where $D$ stands for the Poincar{\'e} duality and $f_*$ for the homological push-forward.

\subsection{The Fock space}

Let $S$ be a smooth complex projective surface with \emph{canonical class}
$K_S$ and \emph{Euler class} $e_S$. For any integer $n\geq 0$, define the
\emph{Hilbert scheme of points} $S^{[n]}:=\Hilb^n(S)$ as the scheme
representing the functor of flat families of length $n$ zero-dimensional closed
subschemes on $S$. By a result of Grothendieck (\cite{G}), it has a natural
structure of projective scheme and is equipped with a \emph{universal family}
$\Xi^S_n\subset S^{[n]}\times S$. By a theorem of Fogarty (\cite{Fo}),
$S^{[n]}$ is in fact a smooth manifold of complex dimension $2n$. We study his
singular cohomology with rational coefficients:
$$
\IH^S_n:=\bigoplus_{i=0}^{4n} H^i(S^{[n]}), \quad \IH^S:=\bigoplus_{n\geq 0}
\IH^S_n.
$$

The unit in $\IH^S_0\cong \IQ$ is called \emph{vaccum vector} and denoted by
$\vac$ (or $\vac_S$ if necessary). The \emph{Fock space} $\IH^S$ is double
graded by $(n,i)$: the integer $n$ is called the \emph{conformal weight} and
the integer $i$ the \emph{cohomological degree}, also denoted by $|\cdot|$.

A linear operator $\kf\in \End(\IH^S)$ is \emph{homogeneous} of bidegree
$(u,v)$ if for any $n$ we have $\kf\left(H^i(S^{[n]})\right)\subset~
H^{i+v}(S^{[n+u]})$. The \emph{super-commutator} of two homogeneous operators
$\kf,\kg$ is defined by:
$$
[\kf,\kg]:=\kf\circ\kg-(-1)^{|\kf|\cdot|\kg|}\kg\circ \kf.
$$

The intersection pairing $\langle\alpha,\beta\rangle_n:=\int_{S^{[n]}}
\alpha\cdot\beta$ for $\alpha,\beta\in \IH^S_n$ extends naturally to a
non-degenerate anti-symmetric bilinear form $\langle\cdot,\cdot\rangle$ on
whole $\IH^S$. For any homogeneous operator $\kf\in \End(\IH^S)$, his
\emph{adjoint} $\kf^\dagger$ is characterized by the relation:
$$
\langle
\kf(\alpha),\beta\rangle=(-1)^{|\kf|\cdot|\alpha|}\langle\alpha,\kf^\dagger(\beta)\rangle.
$$

For $n\geq 0$ and $k>0$, let $S^{[n,n+k]}\subset S^{[n]}\times S\times
S^{[n+k]}$ be the subvariety defined set-theoretically by:
$$
S^{[n,n+k]}:=\left\{(\xi,x,\xi')\,|\,\xi\subset \xi' \text{ and }
\Supp(\cI_{\xi}/\cI_{\xi'})={x}\right\},
$$
where $\cI_\xi$ denotes the ideal sheaf of the subscheme $\xi$ (with
$S^{[n,n]}=\emptyset$). We denote the projections by:
$$
\xymatrix{& S^{[n]}\times S\times
S^{[n+k]}\ar[dl]_\varphi\ar[d]^\rho\ar[dr]^\psi\\S^{[n]}& S& S^{[n+k]}}
$$
The \emph{Heisenberg operators} are the linear operators
$$
\kq_k:H^*(S)\rightarrow \End(\IH^S), \quad k\in \IZ
$$
defined as follows. If $k\geq 0$, for $\alpha\in H^*(S)$ and $x\in H^*(S^{[n]})$ we
set
$$
\kq_k(\alpha)(x):=\psi_!\left(\left[S^{[n,n+k]}\right] \cdot
\varphi^*(x)\cdot \rho^*(\alpha)\right)
$$
and the operators for negative indices are defined by adjonction:
$$
\kq_{-k}(\alpha):=(-1)^k \kq_k(\alpha)^\dagger,\quad \forall k\geq 0.
$$

By convention, $\kq_0=0$. The operators $\kq_k$ are called \emph{creation
operators} if $k\geq 1$ and \emph{annihilation operators} if $k\leq -1$.

\begin{theorem}[Nakajima]\emph{(\cite{N1,N3})} The operators $\kq$ satisfy the
following commutation formula:
$$
[\kq_i(\alpha),\kq_j(\beta)]=i\cdot\delta_{i+j,0}\cdot \int_S\alpha\beta\cdot
\id_{\IH^S}.
$$
\end{theorem}

In particular, the total cohomology space $\IH^S$ admits a basis of vectors
$$
\kq_{n_1}(u^S_1)\cdots\kq_{n_k}(u^S_k)\vac
$$
for $n_i\geq 1$, where the classes $u^S_i$ run over a basis of $H^*(S)$.

For $k\geq 0$, we denote by $\tau_{k!}:H^*(S)\rightarrow H^*(S^k)$ the
push-forward map induced by the diagonal inclusion $\tau_k:S\rightarrow S^k$.
For $k=0$, $\tau_{0!}\alpha$ is understood to be $\int_X \alpha$.
For $k\geq 1$, by K{\"u}nneth decomposition we can set:
\begin{equation}
\label{eq:Kunneth}
\tau_{k!}\alpha=\sum_i \alpha_{i,1}\otimes \cdots\otimes \alpha_{i,k} \in
H^*(S)\otimes \dots \otimes H^*(S).
\end{equation}
We shall make use of the following technical formulas, of the same spirit as in \cite[Lemma $3.1$]{LQW2}:
\begin{lemma} 
\label{lemm:diagonal}
Let $\alpha,\beta,\gamma\in H^*(S)$ and $p,q\geq 1$. Assume that $\tau_{2!}\gamma=\sum\limits_i \gamma_{i,1}\otimes \gamma_{i,2}$. Then:
\begin{align*}
(a)\qquad & \tau_{(k-1)!}\alpha=\sum_i \int_S\alpha_{i,k}\cdot \alpha_{i,1}\otimes \cdots\otimes \alpha_{i,k-1} \\
(b)\qquad & \tau_{(p+q)!}(\alpha\beta\gamma)=\sum_i \tau_{p!}(\alpha \gamma_{i,1})\otimes \tau_{q!}(\beta \gamma_{i,2})
\end{align*}
\end{lemma}

\begin{proof}  The proof follows \cite[Lemma $3.1$]{LQW2}. We recall it for the reader's convenience. For formula (a), consider the following commutative diagram:
$$
\xymatrix{S \ar[r]^-{\tau_k}\ar[dr]^{\tau_{k-1}}& S\times \cdots\times S\ar[d]^{p_{[1,k-1]}}\ar[r]^-{p_k} & S\ar[d]^{\tau_0}\\ & S\times \cdots\times S \ar[r]^-{pr}& {*}}
$$
Then by projection formula one gets:
\begin{align*}
\tau_{(k-1)!}\alpha&=\sum_i p_{[1,k-1]!} \left( \alpha_{i,1}\otimes \cdots\otimes \alpha_{i,k} \right) \\
&= \sum_ip_{[1,k-1]!} \left(p_{[1,k-1]}^* (\alpha_{i,1}\otimes \cdots\otimes \alpha_{i,k-1})\cdot p_k^*\alpha_{i,k} \right) \\
&= \sum_i\alpha_{i,1}\otimes \cdots\otimes \alpha_{i,k-1} \cdot p_{[1,k-1]!}p_k^*\alpha_{i,k} \\
&= \sum_i\alpha_{i,1}\otimes \cdots\otimes \alpha_{i,k-1} \cdot pr^*\tau_{0!}\alpha_{i,k} \\
&=\sum_i \int_S\alpha_{i,k}\cdot \alpha_{i,1}\otimes \cdots\otimes \alpha_{i,k-1}
\end{align*}
For formula (b), notice that:
\begin{align*}
\tau_{2!}(\alpha\beta\gamma)&=\tau_{2!}(\tau_2^*(\alpha\otimes\beta)\cdot\gamma) \\
&=  (\alpha\otimes\beta) \cdot\tau_{2!}(\gamma)\\
&= \sum_i (\alpha\gamma_{i,1})\otimes (\beta\gamma_{i,2})
\end{align*}
and then since $\tau_{p+q}=(\tau_p\times \tau_q)\circ \tau_2$ one gets the result. 
\end{proof}

For $k\geq 1$ and $\alpha\in H^*(S)$, by use of the K{\"u}nneth decomposition (\ref{eq:Kunneth}) we define the \emph{elementary operators} as the operators:
$$
\kq_{n_1}\cdots\kq_{n_k}(\tau_{k!}\alpha):=\sum_i
\kq_{n_1}(\alpha_{i,1})\circ\cdots\circ \kq_{n_k}(\alpha_{i,k}).
$$
The \emph{normally ordered product} of two operators $\kq$ is defined by the
convention:
$$
:\kq_n\kq_m:\;:=\left\{\begin{array}{ll} \kq_n\kq_m&\text{if }n\geq m\\
\kq_m\kq_n&\text{if }n\leq m \end{array}\right.
$$
The \emph{Virasoro operators} are the linear operators
$$
\kL_n:H^*(S)\rightarrow\End(\IH^S), \quad n\in \IZ
$$
defined by $\kL_n:=\frac{1}{2}\sum\limits_{\nu\in \IZ} :\kq_\nu\kq_{n-\nu}:\tau_{2!}$.

\begin{theorem}[Lehn]\emph{(\cite[Theorem 3.3]{L})} The operators $\kL$ satisfy the
following commutation formulas:
\begin{align*}
[\kL_n(\alpha),\kq_m(\beta)]&=-m\cdot\kq_{n+m}(\alpha\beta);\\
[\kL_n(\alpha),\kL_m(\beta)]&=(n-m)\cdot
\kL_{n+m}(\alpha\beta)-\frac{n^3-n}{12}\delta_{n+m,0}\cdot
\int_Se_S\alpha\beta\cdot\id_{\IH^S}.
\end{align*}
\end{theorem}

Denote the canonical projection of the universal family on the Hilbert scheme by $p:\Xi_n^S\rightarrow S^{[n]}$ and let
$B^S_n:=p_*\cO_{\Xi_n^S}$ be the rank $n$ \emph{tautological bundle} on
$S^{[n]}$. Let $\kd\in \End(\IH^S)$ be the
linear operator defined by:
$$
\kd(x):=c_1(B_n^S)\cdot x\quad \forall x\in H^*(S^{[n]}).
$$
The \emph{derivative} of a linear operator $\kf\in \End(\IH^S)$ is defined by
$\kf':=[\kd,\kf]$ and the higher derivatives are $\kf^{(n)}:=(\ad \kd)^n(\kf)$.

\begin{theorem}[Lehn]\emph{(\cite[Main Theorem 3.10]{L})} The derivatives of the operators $\kq$
satisfy the formulas:
\begin{align*}
[\kq'_n(\alpha),\kq_m(\beta)]&=-nm\cdot\left(\kq_{n+m}(\alpha\beta)+\frac{|n|-1}{2}
\delta_{n+m,0}\cdot\int_S K_S\alpha\beta\cdot\id_{\IH^S}\right);\\
\kq'_n(\alpha)&=n\cdot\kL_n(\alpha)+\kq_n(K_S\alpha).
\end{align*}
\end{theorem}

\subsection{Tautological bundles}

Consider the following diagram:
$$
\xymatrix{\Xi_n^S\ar@{^(->}[r]&S^{[n]}\times S\ar[d]_p \ar[r]^-q& S \\&
S^{[n]}}
$$
Let $F$ be a locally free sheaf on $S$. For any $n\geq 0$, the associated
\emph{tautological bundle} on $S^{[n]}$ is defined as:
$$
F^{[n]}:=p_*\left(\cO_{\Xi_n^S}\otimes q^*F\right).
$$
Since the projection $p$ is flat and finite of degree $n$, $F^{[n]}$ is a fibre
bundle of rank $n\cdot \rk(F)$ (by convention, $F^{[0]}=0$). This construction
extends naturally to a well-defined group homomorphism:
$$
-^{[n]}:K(S)\rightarrow K(S^{[n]}).
$$

For $u\in K(S)$, let $\kc(u)$ and
$\kch(u)$ be the linear operators acting for any $n\geq 0$ on $H^*(S^{[n]})$ by
multiplication by the total Chern class $c(u^{[n]})$ and the total Chern
character $ch(u^{[n]})$ respectively.

\begin{theorem}[Lehn]\emph{(\cite[Theorem 4.2]{L})}\label{th:CommCh} Let $u\in K(S)$ be the class of a vector
bundle of rank $r$ and $\alpha\in H^*(S)$. Then:
\begin{align*}
[\kch(u),\kq_1(\alpha)]&=\exp(\ad \kd)(\kq_1(ch(u)\alpha));\\
\kc(u)\circ\kq_1(\alpha)\circ\kc(u)^{-1}&=\sum_{\nu,k\geq
0}\binom{r-k}{\nu}\kq_1^{(\nu)}(c_k(u)\alpha).
\end{align*}
\end{theorem}

By analogy with the construction of tautological bundles, one defines a linear
operation $-^{[n]}:H^*(S)\rightarrow H^*(S^{[n]})$ as follows. For any cohomology class
$\gamma\in H^*(S)$ we set
$$
\gamma^{[n]}:=p_*\left(ch(\cO_{\Xi_n^S})\cdot q^*td(S)\cdot q^*\gamma\right)
$$
where $td(S)$ denotes the Todd class of the tangent bundle $TS$ and we define a linear operator $\kG(\gamma)\in \End(\IH^S)$ acting on
$H^*(S^{[n]})$ by multiplication by $\gamma^{[n]}$.
This definition is such that, by the Riemann-Roch-Grothendieck theorem, the following diagram is commutative:
$$
\xymatrix{H^*(S)\ar[r]^-{-^{[n]}}& H^*(S^{[n]})\\
K(S)\ar[u]^{ch} \ar[r]^-{-^{[n]}}& K(S^{[n]})\ar[u]^{ch}}
$$

\begin{theorem}[Li-Qin-Wang] \emph{(\cite[Lemma 5.8]{LQW3})}
Let $\gamma,\alpha\in~ H^*(S)$. Then:
$$
[\kG(\gamma),\kq_1(\alpha)]=\exp(\ad \kd)(\kq_1(\gamma\alpha)).
$$
\end{theorem}

\section{Universal formulas}

In this section, we develop the notion of \emph{universal formula} and we prove
some general results about the existence and the structure of universal
formulas for the characteristic classes of natural bundles on Hilbert schemes
of points. 

\subsection{Definition of a universal formula}\label{ss:defuniversel}

For any projective variety $X$, let $U^X$ be a cohomology class in $H^*(X)$
which is functorial with respect to pull-backs. For any smooth projective
surface $S$, we set $\cU^S_n:=U^{S^{[n]}}$ and $\cU^S:=\sum\limits_{n\geq 0}
\cU^S_n \in \IH^S$.

\begin{definition}
A class $\cU^S_n\in H^*(S^{[n]})$ admits a \emph{universal formula} if there
exists a polynomial $P\in \IQ[Z_1,\ldots,Z_p]$ independent of $S$, integers
$k_1,\ldots,k_p\geq 1$, indices $n_{i,j}\geq 1$ for $1\leq i \leq p$ and $1\leq
j\leq k_i$ together with cohomology classes $u^S_i\in~ H^*(S)$ in the
sub-algebra generated by $1_S,K_S,e_S$ (which could also depend on some
additional data constructing $\cU^S_n$) such that one has:
$$
\cU^S_n=P\left(\kq_{n_{1,1}}\cdots\kq_{n_{1,k_1}}(\tau_{k_1!}u^S_1),\ldots,
\kq_{n_{p,1}}\cdots\kq_{n_{p,k_p}}(\tau_{k_p!}u^S_p)\right)\vac.
$$
If $P$ is homogeneous of degree $1$, we shall say that the universal formula for
$\cU^S_n$ is a \emph{universal linear combination}.
\end{definition}

Our definition is inspired by the notion of
\emph{universal linear combination} defined by Li-Qin-Wang \cite[Definition
3.1]{LQW2}, but our definition is more restrictive since we only consider
creation operators (for some reasons that will appear soon) and we restrict the
cohomology classes in $H^*(S)$ to the natural sub-algebra generated by the
canonical class and the Euler class.

\subsection{General results on universality}

Let $S_1,S_2$ be two smooth projective surfaces and denote by $S_1\amalg S_2$
their disjoint union. The Hilbert scheme of points decomposes as follows (see
\cite[Formula (0.1)]{EGL}):
$$
(S_1\amalg S_2)^{[n]}=\coprod_{n_1+n_2=n} S_1^{[n_1]}\times S_2^{[n_2]},
$$
inducing the decomposition of the total cohomology:
$$
H^*\left((S_1\amalg S_2)^{[n]}\right)\cong \bigoplus_{n_1+n_2=n}
H^*\left(S_1^{[n_1]}\right)\otimes H^*\left(S_2^{[n_2]}\right),
$$
and the corresponding decomposition of the Fock space:
$$
\IH^{S_1\amalg S_2}\cong \IH^{S_1}\otimes \IH^{S_2}.
$$
In particular, there is a double graded inclusion $\IH^{S_1}\oplus
\IH^{S_2}\subset \IH^{S_1\amalg S_2}$.

\begin{lemma}\label{lemm:univers1}
Suppose that a class $\cU^S\in \IH^S$ admits a universal formula such that:
$$
\cU^{S_1}\oplus \cU^{S_2}=\cU^{S_1\amalg S_2}
$$
for any smooth projective surfaces $S_1,S_2$. Then the universal formula is an
(infinite) universal linear combination.
\end{lemma}

\begin{proof}
Suppose a composition $\kq_{n_1}\cdots \kq_{n_p}(\tau_{p!}u_1^S)\circ
\kq_{m_1}\cdots\kq_{m_q}(\tau_{q!}u_2^S)$ occurs in the universal formula of
$\cU^S$. Then by the decompositions $u_i^{S_1\amalg S_2}=u_i^{S_1}+u_i^{S_2}$
(this is the case for all classes in the sub-algebra generated by
$1_S,K_S,e_S$, and if these classes use additional data, we suppose that this
decomposition holds, as will always be the case in the sequel), the composition
decomposes and extra non-zero terms arise. So the polynomial defining the
universal formula is necessarily of degree $1$.
\end{proof}

\begin{lemma}\label{lemm:univers2}
Suppose that a class $\cU^S\in \IH^S$ admits a universal formula such that
$\cU^S_0=\vac$ and:
$$
\cU^{S_1}\otimes \cU^{S_2}=\cU^{S_1\amalg S_2}
$$
for any smooth projective surfaces $S_1,S_2$. Then the universal formula is an
exponential of an (infinite) universal linear combination.
\end{lemma}

\begin{proof}
Denote by $\kU^S\in \End(\IH^S)$ the operator defined by the universal formula
of $\cU^S$. By construction, $\kU^S\vac=\cU^S$ and with $\vac_{S_1\amalg
S_2}=\vac_{S_1}\otimes \vac_{S_2}$ one gets:
\begin{align*}
\kU^{S_1\amalg S_2}\vac_{S_1\amalg
S_2}&=\left(\kU^{S_1}\vac_{S_1}\right)\otimes
\left(\kU^{S_2}\vac_{S_2}\right) \\
&=\left(\kU^{S_1}\otimes \kU^{S_2}\right)\left(\vac_{S_1}\otimes
\vac_{S_2}\right).
\end{align*}
Since the operator $\kU^S$ contains only creation operators, this equation
implies the equality of the operators:
$$
\kU^{S_1\amalg S_2}=\kU^{S_1}\otimes \kU^{S_2},
$$
and since $\kU^S=\id_{\IH^S}+\cdots$, it admits a logarithm and:
$$
\log \kU^{S_1\amalg S_2}=\log \kU^{S_1} \oplus \log \kU^{S_2}.
$$
Applying the lemma \ref{lemm:univers1} to $\log \kU^S\vac_S$, one gets the
result.
\end{proof}

We dress now a list of some technical results needed to compute with elementary operators; 
the following statements deal with a weaker form of universality (see \S\ref{ss:defuniversel}).

\begin{lemma}[Li-Qin-Wang]\emph{(\cite{LQW2})} \text{}\label{lemm:LQWgen}
\begin{enumerate}
\item \label{lemm:LQWcommute} Any commutator $[\kq_{n_1}\cdots
\kq_{n_p}(\tau_{p!}\alpha),\kq_{m_1}\cdots \kq_{m_q}(\tau_{q!}\beta)]$
can be expressed as a linear combination of elementary operators
$\kq_{i_1}\cdots \kq_{i_k}(\tau_{k!}(\alpha\beta))$, whose coefficients do not
depend on $S$ (with $n_j,m_j,i_j\in \IZ$).

\item \label{lemm:LQWderive} Any derived operator $\kq_n^{(\nu)}(\alpha)$ can be
expressed as a linear combination of elementary operators $\kq_{n_1}\cdots
\kq_{n_k}(\tau_{k!}(K_S^r\alpha))$ for $0\leq r\leq 2$ and $n_i\in \IZ$, whose
coefficients do not depend on $S$.

\item \label{lemm:LQWpassedroite} For any $\alpha\in H^*(S)$, $n_j\in \IZ$, $k\geq 2$
and $1\leq j<k$ one has:
\begin{align*}
&\kq_{n_1}\cdots\kq_{n_j}\kq_{n_{j+1}}\cdots\kq_{n_k}(\tau_{k!}\alpha)
-\kq_{n_1}\cdots\kq_{n_{j+1}}\kq_{n_j}\cdots\kq_{n_k}(\tau_{k!}\alpha)\\
&=n_j\delta_{2n_j+1,0}\kq_{n_1}\cdots\kq_{n_{j-1}}\kq_{n_{j+2}}\cdots\kq_{n_k}(\tau_{(k-2)!}(e_S\alpha))
\end{align*}

\item \label{lemm:LQWGq} Any commutator $[\kG(\alpha),\kq_{n_1}\cdots \kq_{n_k}(\tau_{k!}\beta)]$
can be expressed as a linear combination of operators
$\kq_{m_1}\cdots \kq_{m_p}(\tau_{p!}(K_S^r\alpha\beta))$ with $0\leq r\leq 2$,
whose coefficients do not depend on $S$ (with $n_j$,$m_j\in \IZ$).

\item \label{lemm:LQWcup}
The cup-product of two classes $\kq_{n_1}\cdots\kq_{n_p}(\tau_{p!}\alpha)\vac$
and $\kq_{m_1}\cdots\kq_{m_q}(\tau_{q!}\beta)\vac$ can be expressed as a linear
combination of classes $\kq_{i_1}\cdots\kq_{i_k}(\tau_{k!}\gamma)\vac$ where
$\gamma$ depends on $\alpha,\beta,K_S,e_S$, whose coefficients do not depend on
$S$ (with $n_j$,$m_j$,$i_j\in\IN$).
\end{enumerate}
\end{lemma}

\subsection{Universal formulas for tautological bundles}

\begin{proposition} \label{prop:UniverselChTaut} Let $u\in K(S)$. The Chern
characters $ch(u^{[n]})$ enter in a universal generating series of the kind:
$$
\sum_{n\geq 0} ch(u^{[n]})=\exp(\kq_1(1_S))\kF(u)\vac,
$$
where $\kF(u)\vac$ is an (infinite) universal linear combination depending on
$K_S,e_S,ch(u)$.
\end{proposition}

\begin{proof} Start from the commutation formula given by the theorem
\ref{th:CommCh}:
$$
[\kch(u),\kq_1(1_S)]=\sum_{\nu \geq 0}\frac{1}{\nu !}\kq^{(\nu)}_1(ch(u)),
$$
evaluated in $1_{S^{[n-1]}}$ for $n\geq 1$. Since
$1_{S^{[n]}}=\frac{1}{n!}\kq_1(1_S)^n\vac$, we get:
$$
n\cdot
ch\left(u^{[n]}\right)=\kq_1(1_S)ch\left(u^{[n-1]}\right)+\sum\limits_{\nu \geq
0}\frac{1}{\nu !(n-1)!}\kq^{(\nu)}_1(ch(u))\kq_1(1_S)^{n-1} \vac.
$$
Set $F(t):=\sum\limits_{n\geq 1} ch(u^{[n]})t^n$ (the sum begins in $n=1$ since
$u^{[0]}=0$). Summing up the preceding formula we get:
$$
F'(t)-\kq_1(1_S)F(t)=\left(\sum_{\nu\geq
0}\frac{1}{\nu!}\kq_1^{(\nu)}(ch(u))\right)\exp(\kq_1(1_S)t)\vac.
$$
We have to show that the exponential can be pushed to the left of the formula,
in such a way that the reminding operators form a linear combination of
elementary operators that, applied on the vaccum, can be simplified to a
universal linear combination. This last step is performed with the lemma
\ref{lemm:LQWgen}(\ref{lemm:LQWpassedroite}) which explains how one can push to the right all
annihilation operators occurring in an elementary operator. Then, any annihilation
operator vanishes on the vaccum, so one gets a universal linear combination
(with only creation operators). So it only remains to understand how the
exponential can be pushed to the left.

By lemma \ref{lemm:LQWgen}(\ref{lemm:LQWderive}), any derived operator $\kq_1^{(\nu)}(-)$ is a
linear combination of elementary operators
$\kq_{n_1}\cdots\kq_{n_k}(\tau_{k!}\alpha)$ with coefficients independent of
$S$ so it is enough to study the situation:
$$
\kq_{n_1}\cdots\kq_{n_k}(\tau_{k!}\alpha)\exp(\kq_1(1_S)t)\vac.
$$
If no index $n_i$ is equal to $-1$, there is no problem since the exponential
commutes with $\kq_{n_1}\cdots\kq_{n_k}(\tau_{k!}\alpha)$. Else, by use of the
lemma \ref{lemm:LQWgen}(\ref{lemm:LQWderive}) we can assume that $n_k=-1$ and we observe the
following lemma:

\begin{lemma} \label{lemm:technique} For $\alpha\in H^*(S)$, $x\in \IH^S$
and $n\geq 0$ one has:
\begin{align*}
\kq_{-1}(\alpha)\kq_1(1_S)^nx&=-n\int_S\alpha\cdot
\kq_1(1_S)^{n-1}x+\kq_1(1_S)^n\kq_{-1}(\alpha)x; \\
\kq_{-1}(\alpha)\exp(\kq_1(1_S))x&=-\int_S\alpha\cdot
\exp(\kq_1(1_S))x+\exp(\kq_1(1_S))\kq_{-1}(\alpha)x.
\end{align*}
\end{lemma}

\begin{proof}[Proof of the lemma] The first assertion results on an easy
induction. The second is then straightforward.
\end{proof}

Set the decomposition $\tau_{k!}\alpha=\sum\limits_i \alpha_{i,1}\otimes \cdots \otimes \alpha_{i,k}$.
Then by the preceding lemma applied to $x=\vac$:
\begin{align*}
\kq_{n_1}\cdots \kq_{n_k}(\tau_{k!}\alpha)\circ\exp(\kq_1(1_S))\vac =\sum_i
\kq_{n_1}(\alpha_{i,1})\circ\cdots\circ \kq_{n_k}(\alpha_{i,k})\circ
\exp(\kq_1(1_S))&\vac \\
=-\sum_i \int_S \alpha_{i,k} \cdot\kq_{n_1}(\alpha_{i,1})\circ\cdots\circ
\kq_{n_{k-1}}(\alpha_{i,k-1})\circ \exp(\kq_1(1_S))&\vac.
\end{align*}
Use now the formula (a) of lemma \ref{lemm:diagonal}:
$$
\tau_{(k-1)!}\alpha=\sum_i \int_S\alpha_{i,k} \cdot \alpha_{i,1}\otimes \cdots
\otimes \alpha_{i,k-1},
$$
to get that for $n_k=-1$ one has:
$$
\kq_{n_1}\cdots
\kq_{n_k}(\tau_{k!}\alpha)\circ\exp(\kq_1(1_S))\vac=-\kq_{n_1}\cdots
\kq_{n_{k-1}}(\tau_{(k-1)!}\alpha)\circ\exp(\kq_1(1_S))\vac.
$$
Repeating this process if necessary, one can push the exponential to the left.

We are eventually lead to the following differential equation:
$$
F'(t)-\kq_1(1_S)F(t)=\exp(\kq_1(1_S)t)\kF(u)\vac
$$
where $\kF(u)\vac$ is an universal linear combination. Resolving this equation
we get:
$$
F(t)=\exp(\kq_1(1_S)t)\kF(u)t\vac,
$$
hence the proposition.
\end{proof}

\begin{proposition} \label{prop:UniverselCTaut}
Let $u\in K(S)$. The Chern classes $c(u^{[n]})$ enter in a universal generating
series of the kind:
$$
\sum_{n\geq 0} c(u^{[n]})=\exp(\kF(u))\vac,
$$
where $\kF(u)\vac$ is an (infinite) universal linear combination depending on
$K_S,e_S,c(u)$.
\end{proposition}

\begin{proof}
The Chern classes are polynomials in the Chern characters. By the proposition
\ref{prop:UniverselChTaut}, and the lemma \ref{lemm:LQWgen}(\ref{lemm:LQWcup}) which says that
cup-products of universal formulas are universal, we see that there exists a
universal formula for the Chern classes. Denote by $\cU^S_n:=c(u^{[n]})$ this
universal formula. Let $S_1,S_2$ be two smooth projective surfaces, $n_1,n_2$
two integers and $pr_i:S_1^{[n_1]}\times S_2^{[n_2]}\rightarrow S_i^{[n_i]}$
the projections. For $u_1\oplus u_2\in K(S_1\amalg S_2)=K(S_1)\oplus K(S_2)$
one has the following decomposition (see \cite[Theorem
4.2]{EGL}\footnote{There is an inaccuracy in this formula: instead of a
product $\cdot$ one should read a sum $\oplus$.}):
$$
\left. (u_1\oplus u_2)^{[n_1+n_2]}\right|_{S_1^{[n_1]}\times
S_2^{[n_2]}}=pr_1^*\left(u_1^{[n_1]}\right)\oplus
pr_2^*\left(u_2^{[n_2]}\right),
$$
Hence:
$$
\cU^{S_1\amalg S_2}=\cU^{S_1}\otimes \cU^{S_2}.
$$
Applying the lemma \ref{lemm:univers2}, one gets the result.
\end{proof}

As an example of the formulas we have in mind, let us recall the following
explicit result:
\begin{theorem}[Lehn]\emph{\cite[Theorem 4.6]{L}} Let $L$ be a line bundle
on $S$. Then:
$$
\sum_{n\geq 0}c(L^{[n]})=\exp\left(\sum_{m\geq 1}
\frac{(-1)^{m-1}}{m}\kq_m(c(L))\right)\vac.
$$
\end{theorem}

\subsection{Universal formulas for the tangent bundle}
\label{ss:univTang}

The first result about the universality of the Chern classes of the tangent
bundle $TS^{[n]}$ is the following, where for any partition
$\lambda=(\lambda_1,\ldots,\lambda_k)$ of $2n$, the \emph{Chern numbers} of $S$
are defined as:
$$
c_\lambda(TS^{[n]}):=c_{\lambda_1}(TS^{[n]})\cdots c_{\lambda_k}(TS^{[n]}) \in
H^{4n}(S^{[n]})\cong \IQ.
$$
\begin{theorem}[Ellingsrud-G{\"o}ttsche-Lehn]\emph{\cite[Proposition $0.5$]{EGL}}
For any integer $n$ and any partition $\lambda$ of $2n$, there exists a
universal polynomial $P_\lambda\in \IQ[z_1,z_2]$ such that for any projective
surface $S$ one has:
$$
c_\lambda(TS^{[n]})=P_\lambda\left(c_1(S)^2,c_2(S)\right).
$$
\end{theorem}

We shall prove more precise statements for the Chern classes themselves.

\begin{proposition} \label{prop:UniverselChTan} The Chern
characters $ch(TS^{[n]})$ enter in a universal generating series of the kind:
$$
\sum_{n\geq 0} ch(TS^{[n]})=\exp(\kq_1(1_S))\kF\vac,
$$
where $\kF\vac$ is an (infinite) universal linear combination depending on
$K_S,e_S$.
\end{proposition}

\begin{proof}
Recall some geometric results (see \cite{D,EGL,L}). Consider the following
commutative diagram:
$$
\xymatrix{ & S^{[n,n+1]} \ar[r]^\psi\ar[d]^\sigma \ar@/_/[ddl]_\varphi\ar@/^/[ddr]^\rho & S^{[n+1]} \\
& S^{[n]}\times S \ar[dl]_p \ar[dr]^q\\
S^{[n]}& \Xi_n\ar@{^(->}[u] &S}
$$
Then $\sigma=(\varphi,\rho):S^{[n,n+1]}\rightarrow S^{[n]}\times S$ is the
blow-up of $S^{[n]}\times S$ along the universal family $\Xi_n$, \ie
$S^{[n,n+1]}\cong\Bl_{\Xi_n}(S^{[n]}\times S)$. 
Denote by $E$ the exceptional divisor, set $\cL:=\cO_{S^{[n,n+1]}}(-E)$ and denote the class of the
tangent bundle $TS^{[n]}$ in $K(S^{[n]})$ by $T_n$. Recall the following result:
\begin{proposition}[Ellingsrud-G{\"o}ttsche-Lehn]\emph{\cite[Proposition 2.3]{EGL}}
\label{prop:tangent}
The following relation holds in $K(S^{[n,n+1]})$:
\begin{align*}
\psi^!T_{n+1}=&\varphi^!T_{n}+\cL-\cL\cdot\sigma^!(\cO_{\Xi_n}^\vee)
+\cL^\vee\cdot
\rho^!\omega_S^\vee\\
&-\cL^\vee\cdot \sigma^!(\cO_{\Xi_n})\cdot
\rho^!\omega_S^\vee-\rho^!(\cO_S-T_S+\omega_S^\vee).
\end{align*}
\end{proposition}

Let $\{b_i\}$ be a basis of $H^*(S)$ such that $\int_Sb_ib_jtd(S)=\delta_{i,j}$.
In the K{\"u}nneth decomposition $H^*(S^{[n]}\times S)\cong H^*(S^{[n]})\otimes
H^*(S)$ we can write:
$$
ch(\cO_{\Xi_n})=\sum_i \alpha_i\otimes b_i=\sum_i p^*\alpha_i\cdot q^*b_i
$$
for some classes $\alpha_i$. By definition of the tautological classes and by
use of the projection formula we find:
\begin{align*}
b_j^{[n]}&=p_!\left(ch(\cO_{\Xi_n})\cdot q^*b_j \cdot q^*td(S)\right) \\
&=\sum_i p_!\left(p^*\alpha_i\cdot q^*b_i\cdot q^*b_j \cdot q^*td(S)\right) \\
&=\sum_i \alpha_i\cdot p_!q^*(b_ib_jtd(S) )\\
&=\alpha_j,
\end{align*}
Hence:
$$
ch(\cO_{\Xi_n})=\sum_i b_i^{[n]}\otimes b_i=\sum_i p^*b_i^{[n]}\cdot q^*b_i.
$$
Define ``dual'' tautological classes by:
$$
\gamma^{\{n\}}:=p_!\left(ch(\cO_{\Xi_n}^\vee)\cdot q^*\gamma\cdot
q^*td(S)\right).
$$
One finds similarly:
$$
ch(\cO_{\Xi_n}^\vee)=\sum_i b_i^{\{n\}}\otimes b_i=\sum_i p^*b_i^{\{n\}}\cdot
q^*b_i.
$$

Denote by $\kch T\in \End(\IH^S)$ the operator acting by multiplication by
$ch(T_n)$ on each conformal weight $n$, and by $\kG^\vee(\gamma)$ the operator
multiplying by $\gamma^{\{n\}}$. Note the formulas:
\begin{align*}
\sigma^*ch(\cO_{\Xi_n})&=\sum_i \varphi^*b_i^{[n]}\cdot \rho^*b_i \\
\sigma^*ch(\cO_{\Xi_n}^\vee)&=\sum_i \varphi^*b_i^{\{n\}}\cdot \rho^*b_i,
\end{align*}

We prove the following commutation relation:
\begin{lemma}
\begin{align*}
[\kch T,\kq_1(\alpha)]=&\sum_\nu \frac{1}{\nu!} \kq_1^{(\nu)}(\alpha) \\
&-\sum_{i,\nu} \frac{1}{\nu!} \kq_1^{(\nu)}(b_i\alpha)\circ \kG^\vee(b_i) \\
&+\sum_\nu \frac{(-1)^\nu}{\nu!} \kq_1^{(\nu)}(ch(\omega_S^\vee) \alpha) \\
&-\sum_{i,\nu} \frac{(-1)^\nu}{\nu!} \kq_1^{(\nu)}(b_i
ch(\omega_S^\vee)\alpha)\circ
\kG(b_i) \\
&-\kq_1(ch(\cO_S-T_S+\omega_S^\vee)\alpha).
\end{align*}
\end{lemma}

\begin{proof}[Proof of the lemma] The computation is similar to \cite[Theorem 4.2]{L}, with more terms. 
For any $x\in H^*(S^{[n]})$, by the projection formula one gets:
\begin{align*}
\kch T\circ \kq_1(\alpha)(x)&=ch(T_{n+1})\cdot \psi_!\left( [S^{[n,n+1]}]\cdot \varphi^*(x)\cdot \rho^*(\alpha) \right)\\
&=\psi_!\left( [S^{[n,n+1]}]\cdot \psi^*(ch(T_{n+1}))\cdot\varphi^*(x)\cdot \rho^*(\alpha) \right).
\end{align*}
The proposition \ref{prop:tangent} gives then:
\begin{align*}
\kch T\circ \kq_1(\alpha)(x)=&\psi_!\left( [S^{[n,n+1]}]\cdot\varphi^*(ch(T_{n})\cdot x)\cdot \rho^*(\alpha) \right)\\
&+\psi_!\left( [S^{[n,n+1]}]\cdot ch(\cL) \cdot\varphi^*(x)\cdot \rho^*(\alpha) \right)\\
&-\psi_!\left( [S^{[n,n+1]}]\cdot ch(\cL)\cdot \sigma^*ch(\cO_{\Xi_n}^\vee)\cdot\varphi^*(x)\cdot \rho^*(\alpha) \right)\\
&+\psi_!\left( [S^{[n,n+1]}]\cdot ch(\cL^\vee) \cdot\varphi^*(x)\cdot \rho^*(ch(\omega_S^\vee)\cdot\alpha) \right)\\
&-\psi_!\left( [S^{[n,n+1]}]\cdot ch(\cL^\vee)\cdot \sigma^*ch(\cO_{\Xi_n})\cdot\varphi^*(x)\cdot \rho^*(ch(\omega_S^\vee)\cdot\alpha) \right)\\
&-\psi_!\left( [S^{[n,n+1]}]\cdot\varphi^*(x)\cdot \rho^*(ch(\cO_S-T_S+\omega_S^\vee)\alpha) \right)
\end{align*}
Set $\lambda:=c_1(\cL)$. Then $ch(\cL)=\sum\limits_{\nu\geq 0} \frac{1}{\nu!}\lambda^\nu$ and using the preceding decompositions we find:
\begin{align*}
\kch T\circ \kq_1(\alpha)(x)=&\kq_1(\alpha)\circ \kch T (x)\\
&+\sum_{\nu}\frac{1}{\nu!}\psi_!\left( [S^{[n,n+1]}]\cdot \lambda^\nu \cdot \varphi^*(x)\cdot \rho^*(\alpha) \right)\\
&-\sum_{i,\nu}\frac{1}{\nu!}\psi_!\left( [S^{[n,n+1]}]\cdot \lambda^\nu \cdot \varphi^*(b_i^{\{n\}}\cdot x)\cdot \rho^*(b_i\alpha) \right)\\
&+\sum_{\nu}\frac{(-1)^\nu}{\nu!}\psi_!\left( [S^{[n,n+1]}]\cdot \lambda^\nu \cdot \varphi^*(x)\cdot \rho^*(ch(\omega_S^\vee)\cdot\alpha) \right)\\
&-\sum_{i,\nu}\frac{(-1)^\nu}{\nu!}\psi_!\left( [S^{[n,n+1]}]\cdot \lambda^\nu \cdot \varphi^*(b_i^{[n]}\cdot x)\cdot \rho^*(b_ich(\omega_S^\vee)\alpha) \right)\\
&-\psi_!\left( [S^{[n,n+1]}] \cdot \varphi^*(x)\cdot \rho^*(ch(\cO_S-T_S+\omega_S^\vee)\alpha) \right)
\end{align*}
As explained in the proof of \cite[Theorem 4.2]{L} (or \cite[Lemma 3.9]{L}), a cycle $[S^{[n,n+1]}]\cdot \lambda^\nu$ induces the operator $\kq_1^{(\nu)}$ and with our notations we get:
\begin{align*}
\kch T\circ \kq_1(\alpha)(x)=&\kq_1(\alpha)\circ \kch T (x)\\
&+\sum_{\nu}\frac{1}{\nu!}\kq_1^{(\nu)}(\alpha)(x)\\ 
&-\sum_{i,\nu} \frac{1}{\nu!} \kq_1^{(\nu)}(b_i\alpha)\circ \kG^\vee(b_i)(x)\\ 
&+\sum_\nu \frac{(-1)^\nu}{\nu!} \kq_1^{(\nu)}(ch(\omega_S^\vee) \alpha)(x) \\
&-\sum_{i,\nu} \frac{(-1)^\nu}{\nu!} \kq_1^{(\nu)}(b_i
ch(\omega_S^\vee)\alpha)\circ \kG(b_i)(x) \\
&-\kq_1(ch(\cO_S-T_S+\omega_S^\vee)\alpha)(x).
\end{align*}
\end{proof}

Since $c_1(T_S)=-K_S, c_2(T_S)=e_S,c_1(\omega_S)=K_S$, we deduce
$ch(\omega_S)=~\exp(K_S)$ and $ch(\cO_S-T_S+\omega_S^\vee)=e_S+K_S^2$. We
follow now the same argument as for the proposition \ref{prop:UniverselChTaut}.
We have to show that the following classes are universal formulas for which the
exponential can be pushed to the left:
\begin{align}
&\kq_1^{(\nu)}(1_S) \exp(\kq_1(1_S)) \vac, \label{c1} \\
&\sum_i \kq_1^{(\nu)}(b_i) \kG^\vee(b_i)\exp(\kq_1(1_S)) \vac, \label{c2} \\
&\kq_1^{(\nu)}(\exp(-K_S))\exp(\kq_1(1_S)) \vac, \label{c3} \\
&\sum_i   \kq_1^{(\nu)}(b_i\exp(-K_S))\kG(b_i)\exp(\kq_1(1_S)) \vac, \label{c4} \\
&\kq_1(e_S+K_S^2) \exp(\kq_1(1_S)) \vac. \label{c5}
\end{align}
Formulas (\ref{c1}) and (\ref{c3}) have been studied before, and formula
(\ref{c5}) is obvious. We shall study in details the formula (\ref{c4}), and
then explain how one deduces the result for the formula (\ref{c2}).

By lemma \ref{lemm:LQWgen}(\ref{lemm:LQWderive}), we can assume that
$\kq_1^{(\nu)}(b_i\exp(-K_S))$ is only an elementary operator $\kq_{n_1}\cdots
\kq_{n_k}(\tau_{k!}(b_i\alpha))$ where $\alpha$ is a polynomial in $K_S$. Then:
$$
\kq_{n_1}\cdots \kq_{n_k}(\tau_{k!}(b_i\alpha))
\kG(b_i)\exp(\kq_1(1_S))\vac=\kq_{n_1}\cdots \kq_{n_k}(\tau_{k!}(b_i\alpha))
\sum_{n\geq 0} b_i^{[n]}.
$$
Observe that $b_i^{[n]}=ch\left((ch^{-1}b_i)^{[n]}\right)$ so by the
proposition \ref{prop:UniverselChTaut}, $\sum\limits_{n\geq 0} b_i^{[n]}$
admits a universal formula of the kind\footnote{See also
\cite[Theorem 4.1]{LQW2}.} $\exp(\kq_1(1_S))\kF(b_i)\vac$ where $\kF(b_i)\vac$ is a universal linear
combination depending linearly on $b_i$. The lemma \ref{lemm:technique} applied
to $x=\kF(b_i)\vac$ shows, as in the proof of the proposition
\ref{prop:UniverselChTaut}, that:
$$
\kq_{n_1}\cdots \kq_{n_k}(\tau_{k!}(b_i\alpha))\exp(\kq_1(1_S))\kF(b_i)\vac
$$
can be expressed in a similar form with the exponential on the left, followed
by a linear combination of operators:
$$
\kq_{i_1}\cdots \kq_{i_p}(\tau_{p!}(b_i\beta))\kF(b_i).
$$
By the linearity of the universal formula $\kF(b_i)$, we assume that it
consists only on an elementary operator $\kq_{j_1}\cdots
\kq_{j_q}(\tau_{q!}(b_i\gamma))$. In order to get a universal formula, we have
to get rid of the classes $b_i$. To do so, we call back the summation over the
indices $i$:
\begin{align*}
\sum_i\big(\kq_{i_1}\cdots
\kq_{i_p}(\tau_{p!}(b_i\beta))\big)&\big(\kq_{j_1}\cdots
\kq_{j_q}(\tau_{q!}(b_i\gamma))\big) \\
&=\kq_{i_1}\cdots\kq_{i_p}\kq_{j_1}\cdots\kq_{j_q}\left(\sum_i\tau_{p!}(b_i\beta)\otimes
\tau_{q!}(b_i\gamma)\right)
\end{align*}
Since $\tau_{2!}td(S)=\sum\limits_i b_i\otimes b_i$, formula (b) of lemma \ref{lemm:diagonal} gives
$$
\tau_{(p+q)!}(td(S)\beta\gamma)=\sum_i \tau_{p!}(b_i\beta)\otimes
\tau_{q!}(b_i\gamma)
$$
hence:
\begin{align*}
\sum_i\big(\kq_{i_1}\cdots
\kq_{i_p}(\tau_{p!}(b_i\beta))\big)&\big(\kq_{j_1}\cdots
\kq_{j_q}(\tau_{q!}(b_i\gamma))\big) \\
&=\kq_{i_1}\cdots\kq_{i_p}\kq_{j_1}\cdots\kq_{j_q}(\tau_{(p+q)!}(td(S)\beta\gamma)),
\end{align*}
which shows that we have an universal formula.

The case of the formula (\ref{c2}) is similar since $ch(\cO_{\Xi_n})$ and
$ch(\cO_{\Xi_n}^\vee)$ differ only in some signs for some cohomogical degrees,
so the operators $\kG(\gamma)$ et $\kG^\vee(\gamma)$ behave similarly for all
the results we have used: in the universal formulas, only some signs are
different.
\end{proof}

\begin{proposition} \label{prop:UniverselCTan}
The Chern classes $c(TS^{[n]})$ enter in a universal generating series of the
kind:
$$
\sum_{n\geq 0} c(u^{[n]})=\exp(\kF)\vac,
$$
where $\kF\vac$ is an (infinite) universal linear combination depending on
$K_S,e_S$.
\end{proposition}

\begin{proof}
The proof is similar to the proof of the proposition \ref{prop:UniverselCTaut}
since if $S_1,S_2$ are two smooth projective surfaces, we have the following
decomposition:
$$
\left. T(S_1\amalg S_2)^{[n_1+n_2]}\right|_{S_1^{[n_1]}\times
S_2^{[n_2]}}=pr_1^*\left(TS_1^{[n_1]}\right)\oplus
pr_2^*\left(TS_2^{[n_2]}\right).
$$
\end{proof}

\section{Hilbert scheme of points in the affine plane}

\subsection{Cohomology via symmetric functions}

The \emph{ring of symmetric functions} is the polynomial ring
$\Lambda:=~\IQ[p_1,p_2,\cdots]$ in a countably infinite number of
indeterminates. This space is given a double grading by letting $p_i$ have
\emph{conformal weight} $i$ and \emph{cohomological degree} $i-1$. We denote by
$\Lambda^n$ the subspace of polynomials of conformal weight $n$.

The manifold $\hilb$ has no odd cohomology; his even cohomology has no torsion
and is generated by algebraic cycles (see \cite{ES}). The preceding
construction of Heisenberg operators (naturally extended to the quasi-projective setup by use of the Borel-Moore homology, see \cite{N1,N3}) induces a natural isomorphism
$$
\Lambda^n \cong H^*(\hilb).
$$
In the sequel, we shall study Chern classes of vector bundles on $\hilb$. All
the formulas will be written in the space $\Lambda^n$, since this space
provides a powerful tool for various computations.

A \emph{partition} of an integer $n$ is a decreasing sequence $\lambda:=~(\lambda_1,\ldots,\lambda_k)$ of non-negative
integers  such that $\sum\limits_{i=1}^k\lambda_i=n$ (denoted by $\lambda\vdash n$). The $\lambda_i$
are the \emph{parts} of the partition. The number $l(\lambda)$ of non-zero
parts is the \emph{length} of the partition and the sum $|\lambda|$ of the
parts is the \emph{weight}. If a partition $\lambda$ has $\alpha_1$ parts equal
to $1$, $\alpha_2$ parts equal to $2$, $\ldots$ we shall also denote it by
$\lambda:=(1^{\alpha_1},2^{\alpha_2},\ldots)$ and we set
$z_\lambda:=\prod\limits_{r\geq 1}\alpha_r!r^{\alpha_r}$. A natural basis of
$\Lambda^n$ is given by the \emph{Newton functions}
$p_\lambda:=p_{\lambda_1}\cdots p_{\lambda_k}$ for all partitions $\lambda$ of
$n$.

The \emph{Young diagram} of a partition $\lambda$ is defined by:
$$
D(\lambda):=\{(i,j) \in \IN\times\IN \,|\, j<\lambda_{i+1} \}.
$$
In the representation of such a diagram, we follow a matrix convention:
$$
\Yvcentermath1\young(~xhh,~h~,~) 
\begin{array}{cc}
\lambda=(4,3,1) & x=(0,1)\\
|\lambda|=8 & h(x)=4 \\
l(\lambda)=3 &
\end{array}
$$
where for each cell $x\in D(\lambda)$, the \emph{hook length} $h(x)$ at $x$ is
the number of cells on the right and below $x$ (including the cell $x$ itself) and we set
$h(\lambda):=\prod\limits_{x\in D(\lambda)}h(x)$.

Let $\cC(S_n)$ be the $\IQ$-vector space of \emph{class functions} on $S_n$.
Since conjugacy classes in $S_n$ are indexed by partitions, the functions
$\chi_\lambda$ taking the value $1$ on the conjugacy class $\lambda$ and $0$
else form a basis of $\cC(S_n)$. Let $R(S_n)$ be the $\IQ$-vector space of
representations of $S_n$. By associating to each representation of $S_n$ his
\emph{character}, one gets an isomorphism $\chi:R(S_n)\rightarrow \cC(S_n)$.
The \emph{Frobenius morphism} is the isomorphism $\Phi:\cC(S_n)\rightarrow
\Lambda^n$ characterized by $\Phi(\chi_\lambda)=z_\lambda^{-1}p_\lambda$.
Denote by $\chi^\lambda$ the class function such that
$\Phi(\chi^\lambda)=s_\lambda$ and by $\chi^\lambda_\mu$ the value of
$\chi^\lambda$ at the conjugacy class $\mu$. Then the representations of
characters $\chi^\lambda$ are the irreducible representations of $S_n$.

\subsection{Generating formulas for the tautological bundle on $\hilb$}

The following formulas for the total Chern class and the total Chern character
provide basic examples of the results we have proved about the structure of
universal formulas.

The generating formula for the total Chern class of the tautological bundle
$B_n$ on $\hilb$ is well-known:

\begin{proposition}[Lehn]\emph{(\cite[Theorem 4.6]{L})}
The total Chern class of the tautological bundle on $\hilb$ is given in
$\Lambda$ by the following generating formula:
$$
\sum_{n\geq 0} c(B_n)=\exp\left(\sum_{m\geq 1}(-1)^{m-1} \frac{p_m}{m}\right).
$$
\end{proposition}

The generating formula for the total Chern character is an easy consequence of
a well-known result:

\begin{proposition} \label{prop:GenFormChTaut}
The total Chern character of the tautological bundle on $\hilb$ is given in
$\Lambda$ by the following generating formula:
$$
\sum_{n\geq 0} ch(B_n)=e^{p_1} \sum_{k\geq 1}(-1)^{k-1} \frac{p_k}{k!}.
$$
\end{proposition}

\begin{proof}
By Lehn's theorem (see \cite[Theorem 4.1]{LS}), the Chern character of $B_n$ is
given by:
$$
ch(B_n)=\cD\left(\frac{1}{n!} p_1^n\right),
$$
where the operator $\cD$ is defined by:
$$
\cD=\left. \left( -\sum\limits_{r\geq 1}p_r t^r\right) \exp \left(
-\sum\limits_{r\geq 1}r\frac{\partial}{\partial p_r}t^{-r}\right)
\right|_{t^0}.
$$
Developing the expression we get:
\begin{align*}
\cD &=\left. \left( -\sum\limits_{r\geq 1}p_r t^r \right)
\left(1+\sum\limits_{k\geq 1}\frac{(-1)^k}{k!}\left( \sum\limits_{r\geq
1}r\frac{\partial}{\partial p_r}t^{-r}\right)
^{k}\right) \right|_{t^0} \\
&=\textstyle \left. \left(-\sum\limits_{r\geq 1}p_r t^r\right) \left(
1+\sum\limits_{k\geq 1}\frac{(-1)^k}{k!}\sum\limits_{n_1,\ldots ,n_k\geq 1}
n_1\cdots n_k \frac{\partial}{\partial p_{n_1}} \cdots \frac{\partial}{\partial
p_{n_k}}t^{-(n_1+\cdots +n_k)}\right)
\right|_{t^0} \\
&=\sum\limits_{k\geq 1}\frac{(-1)^{k-1}}{k!} \sum\limits_{n_1,\ldots ,n_k \geq
1}n_1\cdots n_kp_{n_1+\cdots +n_k}\frac{\partial}{\partial p_{n_1}}\cdots
\frac{\partial}{\partial p_{n_k}} .
\end{align*}
Hence:
$$
\cD(p_1^n) =\sum_{k=1}^n (-1)^{k-1}\binom{n}{k} p_1^{n-k}p_k ,
$$
so:
$$
ch(B_n)=\sum_{k=1}^n\frac{(-1)^{k-1}}{k!(n-k)!}p_1^{n-k}p_k ,
$$
hence the formula.
\end{proof}

\subsection{Chern classes of linearized bundles}

The torus $T=\IC^*$ acts on $\IC[x,y]$ by $s.x=sx, s.y=s^{-1}y$ for $s\in T$.
This induces a natural action on $\hilb$ with finitely many fixed points
$\xi_\lambda$ parameterized by the partitions $\lambda$ of $n$.

Let $F$ be a $T$-linearized bundle of rank $r$ on $\hilb$. Each fibre
$F(\xi_\lambda)$ has a structure of representation of $T$, uniquely determined
by its \emph{weights} $f_1^\lambda,\ldots,f_r^\lambda$. These data are enough
to recover the Chern classes and the Chern characters of $F$:

\begin{theorem}[Boissi{\`e}re]\emph{\cite[Theorem 4.2]{SB1}} \label{th:421}
Let $F$ be a $T$-linearized vector bundle of rank $r$ on $\hilb$ and
$f_1^\lambda,\ldots,f_r^\lambda$ the weights of the action on the fibre at each
fixed point. Then the Chern classes of $F$ written in $\Lambda^n$ are:
$$
c_k(F)=\sum\limits_{\lambda\vdash n} \frac{1}{h(\lambda)}
\Coeff\left(t^k,\prod_{i=1}^r(1+f_i^\lambda t)\right) \sum\limits_{\substack{\mu\vdash n \\
l(\mu)=n-k}} z_\mu^{-1}\chi^\lambda_\mu p_\mu.
$$

The Chern characters of $F$ are:
$$
ch_k(F)=\frac{1}{k!}\sum\limits_{\lambda\vdash n} \frac{1}{h(\lambda)}
\sum\limits_{i=1}^r \left(f_i^\lambda\right)^k \sum\limits_{\substack{\mu\vdash n \\
l(\mu)=n-k}} z_\mu^{-1}\chi^\lambda_\mu p_\mu.
$$
\end{theorem}

\section{Chern classes of the tangent bundle on $\hilb$}

We can now prove the two formulas announced in the introduction.

\begin{theorem}
The total Chern class of the tangent bundle on $\hilb$ is given in $\Lambda$ by
the following generating formula :
$$
\sum_{n\geq 0} c(T\hilb)=\exp\left(\sum_{k\geq 0}(-1)^kC_k
\frac{p_{2k+1}}{2k+1}\right),
$$
where $C_k:=\displaystyle\frac{1}{k+1}\binom{2k}{k}$ is the $k$-th Catalan
number.
\end{theorem}

\begin{proof} The proof follows the strategy explained in the introduction.

\textsc{Step 1.} The manifold $\hilb$ is not projective, so in order to use our
general results on universal formulas, we embed it in $\mathbb{P}_2^{[n]}$. By
proposition \ref{prop:UniverselCTan}, the Chern classes of the tangent space on
$\mathbb{P}_2^{[n]}$ are given by a generating formula:
$$
\sum_{n\geq 0} c\left(T\mathbb{P}_2^{[n]}\right)=\exp\left(\kF\right)\vac,
$$
where $\kF$ is a linear combination of elementary operators $\kq_{n_1}\cdots
\kq_{n_k}(\tau_{k!}\alpha)$ (with $n_i\geq 1$) whose argument depends on
$1_{\mathbb{P}_2}$, $K_{\mathbb{P}_2}$ and $e_{\mathbb{P}_2}$. The inclusion
$\IC^2\subset \mathbb{P}_2$ induces an open immersion  $\hilb\subset
\mathbb{P}_2^{[n]}$ giving a surjection $\IH^{\mathbb{P}_2}\twoheadrightarrow
\IH^{\IC^2}$. Since the restrictions to the affine plane of the classes
$K_{\mathbb{P}_2}$ and $e_{\mathbb{P}_2}$ are trivial and the morphisms
$\tau_{k!}$ are zero if $k\geq 2$, the formula for the Chern class of the
tangent bundle is more simple and, denoting $\kq_m:=\kq_m(1_{\IC^2})$ and using
the identification between $\kq_m$ and the Newton function $p_m$, we see that
the formula is only:
$$
\sum_{n\geq 0} c(T\hilb)=\exp\left( \sum_{m\geq 1} f_m p_m \right).
$$

\textsc{Step 2.} In the cohomology $H^*(\hilb)$, the greater non-zero
cohomological degree is $n-1$ and is generated (through the identification with
$\Lambda^n$) by the Newton function $p_n$, the only Newton function $p_\lambda$
of conformal weight $n$ and cohomological degree $n-1$. So if we develop the
exponential in the preceding formula and compare both conformal weights and
cohomological degrees, we see that:
$$
\sum_{n\geq 0} c_{n-1}(T\hilb)= \sum_{m\geq 1} f_m p_m.
$$
So we only have to compute explicitly these ``maximal'' Chern classes.

\textsc{Step 3.} In order to compute the Chern classes $c_{n-1}(T\hilb)$, we
use the theorem \ref{th:421}. The weights of the fibre of the tangent bundle on
$\hilb$ at a fixed point $\xi_\lambda$ are the hook lengths and their
opposites: $\{h(x),-h(x)\,|\,x\in D(\lambda)\}$ (see for example \cite{N1}),
hence:
$$
c_{n-1}(T\hilb)=\sum_{\lambda\vdash n}
\frac{1}{h(\lambda)}\Coeff\left(t^{n-1},\prod_{x\in D(\lambda)}
(1-h(x)^2t^2)\right) \sum_{\substack{\mu\vdash n\\l(\mu)=1}}\chi^\lambda_\mu
z_\mu^{-1}p_\mu.
$$
We deduce that if $n$ is even, then $c_{n-1}(T\hilb)=0$. Suppose now that
$n=2k+1$ and set $\alpha_k:=c_{2k}(T\Hilb^{2k+1}(\IC^2))$. The only partition
$\mu$ of $2k+1$ of length $1$ is $\mu=(2k+1)$ and the evaluation of a character
$\chi^\lambda$ on a maximal cycle follows the following rule (see
\cite[Exercise 4.16]{FH}):
$$
\chi^\lambda_{(2k+1)}= \left\{
\begin{array}{ll}
(-1)^s & \text{if } \lambda=(2k+1-s,1,\ldots,1), \quad 0\leq s\leq 2k \\
0& \text{else}
\end{array}
\right.
$$
For such a partition $\lambda=(2k+1-s,1,\ldots,1)$, the hook lengths are the
integers $\{1,\ldots,s\}$, $\{1,\ldots,2k-s\}$ and $2k+1$ so:
\begin{align*}
\textstyle{\alpha_{k}} &\scriptstyle= \sum\limits_{s=0}^{2k}
\frac{(-1)^s}{(2k+1)s!(2k-s)!}
\Coeff\left(t^{2k},\prod\limits_{i=1}^{s}(1-i^2t^2)\prod\limits_{j=1}^{2k-s}(1-j^2t^2)(1-(2k+1)^2t^2)\right)
\frac{p_{2k+1}}{2k+1} \\
&=\textstyle\Coeff\left(t^{2k},\frac{(1-(2k+1)^2t^2)}{2k+1}\sum\limits_{s=0}^{2k}
\frac{(-1)^s}{s!(2k-s)!}
\prod\limits_{i=1}^{s}(1-i^2t^2)\prod\limits_{j=1}^{2k-s}(1-j^2t^2)\right)\frac{p_{2k+1}}{2k+1}.
\end{align*}

Extract from this formula the following polynomial:
$$
P_k:=\sum\limits_{s=0}^{2k} \frac{(-1)^s}{s!(2k-s)!}
\prod\limits_{i=1}^{s}(1-i^2t^2)\prod\limits_{j=1}^{2k-s}(1-j^2t^2).
$$
Introducing the \emph{Pochhammer symbol}:
$$
(a)_r:=a(a+1)\cdots(a+r-1)  \text{ for } r\geq 1
$$
we get a nicer expression for the polynomial $P_k$:
$$
P_k=t^{4k}\sum\limits_{s=0}^{2k} (-1)^s
\frac{\left(1-\frac{1}{t}\right)_s\left(1+\frac{1}{t}\right)_s}{s!}\frac{\left(1-\frac{1}{t}\right)_{2k-s}\left(1+\frac{1}{t}\right)_{2k-s}}{(2k-s)!}
$$
Observe then the following lemma:

\begin{lemma} For any $a,b$ and any integer $k\geq 0$ we have the following identity:
$$
\sum_{s=0}^{2k} (-1)^s
\frac{(a)_s(b)_s}{s!}\frac{(a)_{2k-s}(b)_{2k-s}}{(2k-s)!}
=\frac{1}{k!}\frac{(a+b)_{2k}(a)_k(b)_k}{(a+b)_k.}
$$
\end{lemma}

\begin{proof}[Proof of the lemma]
Recall the definition of a \emph{generalized hypergeometric function}:
$$
_pF_q(a_1,\ldots,a_p;b_1,\ldots,b_q;z):=\sum_{n\geq 0} \frac{(a_1)_n\cdots
(a_p)_n}{(b_1)_n\cdots (b_q)_n}\frac{z^n}{n!}.
$$
Our identity is then nothing else than the developed form of the following
product identity on hypergeometric functions (see \cite[\S II.2.9, p.
63]{MOS}):
$$
_2F_0(a,b;-z)_2F_0(a,b;z)=
{_4F_1}\left(a,b,\frac{a+b}{2},\frac{a+b+1}{2};a+b;4z^2\right).
$$
\end{proof}

Applying this lemma with $a=1-\frac{1}{t}$ and $b=1+\frac{1}{t}$, we get for
the polynomial $P_k$ the simple expression:
$$
P_k=(-1)^k(2k+1)C_k t^{2k} \prod_{i=1}^k (1-i^2t^2),
$$
from which follows immediately:
$$
\alpha_k=(-1)^kC_k \frac{p_{2k+1}}{2k+1},
$$
hence the theorem.
\end{proof}

\begin{remark}
The preceding theorem was first conjectured, thanks to the inspired impulsion
of Manfred Lehn, by some experiments and numerology. We also thank Marc
Nieper-Wi{\ss}kirchen for his ``hypergeometric'' help.
\end{remark}

\section{Chern Character of the tangent bundle on $\hilb$}

\begin{theorem} \label{th:GenFormChTan}
The Chern characters of the tangent bundle on $\hilb$ are given in $\Lambda$ by
the following generating formula:
$$
\sum_{n\geq 0} ch(T\hilb)=2e^{p_1}\sum_{k\geq 0}\frac{p_{2k+1}}{(2k+1)!}.
$$
\end{theorem}

\begin{proof} The proof is very similar to the preceding proof for the Chern
classes, so we only mention the main differences.

\textsc{Step 1.} As before, we see that there exists an universal formula
$\exp(\kq_1(1_S))\kF\vac$ where $\kF$ is a linear combination which, in the
case of the affine plane, reduces to:
$$
\sum_{n\geq 0}ch(T\hilb)=e^{p_1}\sum_{m\geq 1}f_m p_m.
$$

\textsc{Step 2.} Comparing both conformal weights and cohomological degrees we
see that:
$$
\sum_{n\geq 1}ch_{n-1}(T\hilb)=\sum_{m\geq 1}f_m p_m.
$$

\textsc{Step 3.} As before, we compute the Chern characters $ch_{n-1}(T\hilb)$
with the theorem \ref{th:421}:
$$
\textstyle{ch_{n-1}(T\hilb)=\frac{1}{(n-1)!}\sum\limits_{\lambda\vdash n}
\frac{1}{h(\lambda)} \left(\sum\limits_{x\in D(\lambda)}
(1+(-1)^{n-1})h(x)^{n-1}\right) z_{(n)}^{-1}\chi^\lambda_{(n)} p_n}.
$$
We deduce that $ch_{n-1}(T\hilb)=0$ for $n$ even. If $n=2k+1$ is odd, we set
$\beta_k:=ch_{2k}(T\Hilb^{2k+1}(\IC^2))$ and get:
$$
\beta_k=2\sum\limits_{s=0}^{2k}\frac{(-1)^s}{(2k+1)s!(2k-s)!}\left(\sum\limits_{i=1}^s
i^{2k}+\sum\limits_{j=1}^{2k-s}j^{2k}+(2k+1)^{2k}\right)\frac{p_{2k+1}}{(2k+1)!}.
$$
One verifies that:
\begin{align*}
&\sum\limits_{s=0}^{2k}\frac{(-1)^s}{s!(2k-s)!}=0, \\
&\sum\limits_{s=0}^{2k}\frac{(-1)^s}{s!(2k-s)!}\sum\limits_{i=1}^s i^{2k}
=\frac{1}{2k\cdot(2k)!}\sum_{i=0}^{2k} (-1)^i \binom{2k}{i}i^{2k+1}, \\
&\sum\limits_{s=0}^{2k}\frac{(-1)^s}{s!(2k-s)!}\sum\limits_{j=1}^{2k-s} j^{2k}
=\frac{1}{2k\cdot(2k)!}\sum_{j=0}^{2k} (-1)^j \binom{2k}{j}j^{2k+1},
\end{align*}
by use of the elementary formulas $\sum\limits_{k=0}^p
(-1)^k\binom{n}{k}=(-1)^p\binom{n-1}{p}$ and
$\binom{n-1}{p-1}=\frac{p}{n}\binom{n}{p}$.

Observe then the following lemma:

\begin{lemma} For any $n\geq 0$, the following identity holds:
$$
\sum_{j=0}^n (-1)^{n-j}\binom{n}{j}j^{n+1}=\binom{n+1}{2}n!
$$
\end{lemma}

\begin{proof}[Proof of the lemma]
The \emph{Stirling number of the second kind} $S(n,k)$ is the number of
partitions in $k$ blocks of a set of $n$ elements, and is given by the
following formula (see \cite{Com}):
$$
S(n,k)=\frac{1}{k!} \sum_{j=0}^k (-1)^j \binom{k}{j}(k-j)^n.
$$
In particular:
$$
S(n+1,n)=\frac{1}{n!}\sum_{j=0}^n (-1)^{n-j} \binom{n}{j} j^{n+1},
$$
and we observe that $S(n+1,n)=\binom{n+1}{2}$ since putting $n+1$ identical objects in
$n$ boxes is the same as choosing $2$ objects among $n+1$.
\end{proof}

Applying this lemma we get:
\begin{align*}
&\sum\limits_{s=0}^{2k}\frac{(-1)^s}{s!(2k-s)!}\sum\limits_{i=1}^s i^{2k}
=\frac{2k+1}{2}, \\
&\sum\limits_{s=0}^{2k}\frac{(-1)^s}{s!(2k-s)!}\sum\limits_{j=1}^{2k-s} j^{2k}
=\frac{2k+1}{2},
\end{align*}
from which follows immediately:
$$
\beta_k=2\frac{p_{2k+1}}{(2k+1)!},
$$
hence the theorem.
\end{proof}

Comparing the generating formulas for the Chern character of $T\hilb$ (theorem
\ref{th:GenFormChTan}) and of $B_n$ (theorem \ref{prop:GenFormChTaut}), one
remarks the following decomposition:
\begin{corollary} The following identity holds in $K(\hilb)$:
$$
T\hilb=B_n+B_n^*.
$$
\end{corollary}

\begin{remark}
This formula appears as a particular case of the general decomposition
formula of Ellingsrud-G{\"o}ttsche-Lehn for the tangent bundle over $S^{[n]}$ (see the proof of \cite[Proposition 2.2]{EGL}):
$T_n=B_n+B_n^*-p_!(\cO_n^\vee\cdot \cO_n)$.
\end{remark}

\nocite{*}
\bibliographystyle{amsplain}
\bibliography{BiblioTangentBundle}

\end{document}